\documentclass[leqno,12pt]{article}
\usepackage{amssymb}
\usepackage{amsmath}
\usepackage{amsfonts}
\usepackage{amsthm}
\usepackage{amscd}
\theoremstyle{plain}
\numberwithin{equation}{section}
\newtheorem{thm}{Theorem}[section]

\newtheorem{lem}{Lemma}[subsection]
\newtheorem{corollary}[lem]{Corollary}
\newtheorem{prop}[lem]{Proposition}
\newtheorem{obs}[lem]{Observation}

\theoremstyle{definition}
\newtheorem{rem}[lem]{Remark}
\newtheorem{defn}[lem]{Definition}
\newtheorem{exa}[lem]{Example}
\newtheorem{exam}{Example}[section]

\newenvironment{newenumerate}{
 \begin{enumerate}
 \renewcommand{\labelenumi}{\upshape(\roman{enumi})}
 }{\end{enumerate}}
\renewcommand{\labelenumi}{\theenumi)}
\newfont{\bg}{cmr10 scaled\magstep3}
\newcommand{\bigzero}{\smash{\hbox{\bg 0}}}
\newcommand{\R}{\mathbb R}

\newfont{\Bg}{
 cmr10 scaled\magstep5}

\def \Ad{\operatorname{Ad}}

\def \rank{\operatorname{rank}}
\def \Ker{\operatorname{Ker}}
\def \Hom{\operatorname{Hom}}
\def \trans#1{{}^t\!{#1}}
\def \ya#1{\overset\rightarrow{{#1}}}
\def\rarrowsim{\smash{\mathop{\,\longrightarrow\,}\limits
  ^{\lower1.5pt\hbox{$\scriptstyle\sim$}}}} %TK
\newcommand{\tG}{G \times \dots \times G}
\newcommand{\rdash}{(\mathrm{R})'}
\newcommand{\sdash}{(\mathrm{S})'}
\providecommand{\bysame}{\makebox[3em]{\hrulefill}\thinspace}
\newcommand \ro{\operatorname{ro}}

\begin{document}

\title{Deformation of Properly Discontinuous \\ Actions 
of ${\mathbb{Z}}^k$ 
on $\mathbb{R}^{k+1}$ }
\author{Toshiyuki KOBAYASHI\\
Research Institute for Mathematical Sciences,\\
Kyoto University, Kyoto, 606-8502, Japan\\
{\small \textit{E-mail address}: \texttt{toshi@kurims.kyoto-u.ac.jp}}
\\[\medskipamount]
and
\\[\medskipamount]
Salma NASRIN\\
Graduate School of Mathematical Sciences,\\
University of Tokyo, Tokyo 153-8914, Japan\\
{\small \textit{E-mail address}: \texttt{nasrin@ms.u-tokyo.ac.jp}}}
\date{}
\maketitle

\noindent
\textit{%
Mathematics Subject Classifications\/} (2000): 
Primary:
57S30; \\ %(1980-now) Discontinuous groups of transformations
Secondary:  
%22F30, %Homogeneous spaces
%[For general actions on manifolds or preserving geometrical structures, 
%see 57M60, 57Sxx; for discrete subgroups of Lie groups see especially 22E40]  
22E25, % Nilpotent and solvable Lie groups
22E40,  %(1973-now) Discrete subgroups of Lie groups [See also 20Hxx, 32Nxx]
53C30,  %(1973-now) Homogeneous manifolds
58H15  %(1980-now) Deformations of structures

\medskip
\noindent
\textit{%
Key words\/}:  discontinuous group, 
deformation, rigidity, proper action, 
affine transformation, 
properly discontinuous action, 
homogeneous space.

\medskip\medskip

\begin{abstract}
We consider the deformation of a discontinuous group acting on the
Euclidean space by affine transformations.
A distinguished feature here is that even
a `small' deformation of a 
discrete subgroup may destroy 
proper discontinuity of its action.
In order to understand the local structure of the deformation space 
of discontinuous groups, 
we introduce the concepts from a group
theoretic perspective, and
 focus on `stability' and 
`local rigidity' of discontinuous groups.
As a test case,
we give an explicit description
of the deformation 
space of \ $\mathbb{Z}^k$ acting properly discontinuously on 
$\mathbb{R}^{k+1}$
 by affine nilpotent transformations.
Our method uses an idea of `continuous
 analogue' and relies on the criterion
 of proper actions
  on nilmanifolds.
\end{abstract}

\section{Local rigidity and stability}
Our concern in this article is with the deformation of 
 discontinuous
 groups for non-Riemannian homogeneous spaces.

\subsection{Deformation of discontinuous groups
\\ 
---  the non-Riemannian case} \label{SS:11}

In contrast to the traditional case of
 discontinuous groups
 acting on Riemannian manifolds as isometries,
 our problem in the non-Riemannian case
  includes the following subproblem:
 if a discrete subgroup can be deformed,
 determine the range of the deformation
  parameters that does not destroy
  the proper discontinuity of its action.
   
As a clue to understanding the local structure 
 of the `deformation space',
 we consider a manifold $X$ acted on by
 a Lie group $G$.
Suppose $\Gamma$ is a discontinuous group for 
 $X$,
 that is,
 $\Gamma$ is a discrete subgroup of $G$
 acting properly discontinuously and freely on $X$. 
Let $\Gamma'$ be another discrete subgroup of $G$
 which is \lq{sufficiently close}\rq\ to $\Gamma$. 
Now, our basic question is if the following statements hold or not:
\begin{list}{}{\setlength{\labelwidth}{3em}
\setlength{\leftmargin}{3.2em}}
\item[\mathversion{bold}$\rdash$ :]
(\textit{Local Rigidity})
$\Gamma'$ is conjugate to $\Gamma$
 by an inner automorphism of $G$. 
\item[\mathversion{bold}$\sdash$ :]
(\textit{Stability})
The $\Gamma'$-action on $X$ is properly discontinuous and free.  
\end{list}

For a homogeneous space $X = G/H$
  ($H$ being a closed subgroup of $G$),
  obvious remarks are:
\begin{enumerate}
\item[1)]
  If  $H$ is compact, then $\sdash$ automatically holds.
\item[2)]
If $\rdash$ holds, so does $\sdash$. 
\end{enumerate}

This article studies
 the deformation of discontinuous groups
 for $G/H$
  in the case that $\sdash$ does not hold. 
This implies 
particularly that $H$ is non-compact and 
that $\rdash$ does not hold. 

\smallskip
Let us now formalize the above two statements
 $\rdash$ and $\sdash$ more rigorously.
We begin with an (abstract) finitely generated
 group $\Gamma$,
  and denote by $\Hom(\Gamma, G)$ the set of
 all group homomorphisms of $\Gamma$
  into a Lie group $G$.
Taking generators $\gamma_1, \dots, \gamma_k$ of
 $\Gamma$, we use the injective map
$$
\Hom(\Gamma, G) \hookrightarrow
  G \times \dots \times G,
  \quad
 \varphi \mapsto (\varphi(\gamma_1), \dots,
  \varphi(\gamma_k))
$$
 to endow $\Hom(\Gamma, G)$ with the induced
  topology from the
  direct product $G \times \dots \times G$. 
This topology
 is independent of the choice of generators.
Now, suppose $G$ acts continuously on 
 a manifold $X$.
We recall from \cite{kobayashi93}
 the following definition:
\begin{align}
\label{eqn:R}
R(\Gamma, G; X)  := 
  \{ &\varphi \in \Hom(\Gamma, G) \  : \ 
 \text{$\varphi$ is injective, and  }  \\ 
 & \quad \text{$\varphi(\Gamma)$ acts properly discontinuously 
and freely on $X$} \, \}.\kern-1em
\nonumber
\end{align} 
We remark that our notation here is
 slightly different from that in
 \cite{kobayashi93,kobayashi00}:
 for a homogeneous space $X =G/H$
 our notation $R(\Gamma, G; G/H)$ here
 coincides with $R(\Gamma, G, H)$ loc. cit.

For each $\varphi \in R(\Gamma, G; X)$, the quotient space 
$\varphi(\Gamma)\backslash X \simeq \varphi(\Gamma) \backslash G/H$ 
(a \textit{Clifford--Klein form} of $X$) becomes a Hausdorff 
topological space,
 and can be given a unique manifold structure
 for which the natural quotient map 
 $X \rightarrow \varphi(\Gamma)\backslash X$ is a local diffeomorphism. 
Therefore, 
the Clifford--Klein form $\varphi(\Gamma)\backslash X$ 
enjoys all $G$-invariant local geometric structures on $X$. 
Thus, $R(\Gamma, G; X)$ may be regarded as a parameter space of 
Clifford--Klein forms $\varphi(\Gamma)\backslash X$ 
with parameter $\varphi$.

\medskip 
To be more precise on `parameter',
 we should take into account `unessential' deformation arising from
 inner automorphisms.
If two homomorphisms $\varphi_1$ and $\varphi_2$
 belonging to $R(\Gamma, G; X)$ satisfy
 $\varphi_2 = g \circ \varphi_1 \circ g^{-1}$ for some $g \in G$, 
 then the corresponding
  Clifford--Klein forms are isomorphic to each other by 
the natural diffeomorphism $\varphi_{1} (\Gamma) \backslash X \rarrowsim 
\varphi_{2} (\Gamma) \backslash X,$ 
$\varphi_{1} (\Gamma) x H \mapsto \varphi_{2} (\Gamma) g x H. $ 
In light of this observation,
 we define the \textit{deformation space} 
 as the quotient set 
\begin{equation}
\label{eqn:T}
\mathcal{T} (\Gamma, G; X) := R(\Gamma, G; X)/G. 
\end{equation} 
For example, if $G = PSL(2, \mathbb R)$ and $X$ is the upper half 
plane, and if $\Gamma$ is the fundamental group of a closed 
Riemann surface $M_g$ 
of genus $g \ge 2$, then $\mathcal{T}(\Gamma, G; X)$
  is nothing other than the Teichm\"uller space of $M_g$.
We refer the reader to an expository paper 
\cite{kobayashi00}
for some elementary examples of the deformation space for
non-Riemannian $X$.

\smallskip
Suppose now that $\varphi_0: \Gamma \to G$ belongs to $R(\Gamma, G; X)$,
 and we reformalize
  $\rdash$ and $\sdash$ as follows:
\begin{itemize}
\item[\bfseries(R):] (\textit{Local rigidity})
$G \cdot \varphi_0$ is open in $\operatorname{Hom}(\Gamma, G)$.
\item[\bfseries(S):] (\textit{Stability})
There is an open subset $V$ of $\Hom(\Gamma,G)$ such that
$\varphi_0\in V \subset R(\Gamma,G;X)$.
\end{itemize}

We say $\varphi_0 \in R(\Gamma, G; X)$ 
is \textit{locally rigid as a 
discontinuous group} for $X$ if (R) holds. 
For a 
Riemannian symmetric space $X$, our terminology here is 
consistent with Weil's terminology used in \cite{weil}.

A celebrated Selberg--Weil rigidity \cite{weil} 
 for an irreducible {\bf{Riemannian}}
  symmetric space $X$  asserts
 that
 (R) holds
 for any torsion free uniform lattice
   $\varphi_0(\Gamma)$ of $G$
  unless $G$ is locally isomorphic to $SL(2,{\mathbb{R}})$,
 whereas $R(\Gamma,G;X)$ is always
 open in $\operatorname{Hom}(\Gamma,G)$ 
 and thus (S) holds.  
In contrast, 
there is an example that (R) fails for 
 an irreducible {\bf{non-Riemannian}} symmetric space $X$
  of an arbitrarily high dimension (see \cite{kobayashi93}). 

The failure of (R) arouses our interest in
  the deformation space $\mathcal{T} (\Gamma, G; X)$ like
   the classical Teichm\"uller theory of Riemann surfaces.
Besides,
 the concept of the stability (S) may be regarded as a first 
 step to understand the local structure of the 
 deformation space $\mathcal{T} (\Gamma, G; X)$ 
in the setting where (R) fails, 
 in particular, where $H$ is non-compact.

Such a viewpoint traces back to
 the paper \cite{goldman} 
 on three dimensional Lorentz space forms,
where
 Goldman discovered a discontinuous group $\Gamma$ for which (R) fails,
 and raised a question if (S) still holds
  (not exactly in the way formulated 
 here).
His case concerns with a semisimple Lie group
 $G$ which is locally isomorphic to
  $SO(2,1) \times SO(2,1)$.
This question was solved affirmatively,
 namely, there exists a cocompact
 discontinuous group $\Gamma$ 
 for which (R) fails but (S) holds 
 in the generality that $G$ is
  a semisimple Lie group 
 which is locally isomorphic
 to the direct product of two copies
 of $SO(n,1)$ or $SU(n,1)$
  %where $H \approx SO(n,1)$ or $SU(n,1)$
 (see \cite{kobayashi98, salein}). 
Its proof relies
 on the criterion for properly discontinuous actions
    \cite{benoist, kobayashi96}
on homogeneous spaces of reductive groups. 

\smallskip
For a more general $(\Gamma, G, X)$
 such as the affine transformation group $G$,
 both (R) and (S) can fail,
 as is seen by the following one dimensional example:

\begin{exam} %\label{R:affineone}
Let $\Gamma := \mathbb Z$,
 and
$G$ be the $a x + b$ group,
 that is, $G = \{(a,b):a>0, b \in \mathbb R\}$
 with the multiplication given by
  $(a, b) \cdot (a', b') = (a a',
 a b' + b)$.
Consider the affine transformation of $G$ 
on $X := \mathbb R$.
Then, $\Hom(\Gamma, G) \simeq G$,
 whereas $R(\Gamma, G; X) \simeq \{(1,b):b \neq 0\}$.
Hence, neither (R) nor (S) holds.
\end{exam}

The above example deals with a homogeneous
 space $X$ of a solvable group $G$
 and with a cocompact discontinuous group $\Gamma$.

%\medskip

\subsection{Summary of this article} 
\label{SS:12}

This article analyses the failure of 
  (R) and (S) 
  for homogeneous spaces of
   {\bf nilpotent} Lie groups $G$.
For a simply-connected nilpotent Lie group $G$,
 $R(\Gamma, G; X)$ becomes open in $\Hom(\Gamma, G)$
  and thus (S) always holds
  if $\Gamma \backslash X$ is compact
 (\cite{yoshino}).
Thus, our interest here is in the case when 
 $\Gamma \backslash X$ is non-compact.
As a test case, 
 we initiate a detailed analysis
 on the deformation space
 in the following  setting:
\begin{align*}
&&
\Gamma &:= \mathbb{Z}^k  &&\text{(free abelian group of rank $k$)},
&&
\\
&&
X &:= \mathbb{R}^{k+1}   &&\text{(nilmanifold)}, 
\\
&&
G & \subset \operatorname{Aff}(\mathbb R^{k+1})
 &&\text{(a two-step nilpotent subgroup)},
\end{align*}
where $\Gamma$ acts on $X$ as nilpotent affine
 transformations via $G$.
Then, we propose a method of giving
 a concrete description of
 $R(\Gamma, G; X)$ 
 and the deformation space 
 $\mathcal{T}(\Gamma, G; X)$
 for a specific choice of $G$.
 %(see Proposition \ref{P:homomorphism}, 
 Our main results are 
Theorems \ref{T:deformpara} and \ref{T:deformspace}.

Besides, we shall see in Corollary~\ref{cor:dimdefo} that 
the deformation space contains a smooth manifold
$\mathcal{T}'(\Gamma,G;X)$
as its open dense subset such that
$$
\dim \mathcal{T}'(\Gamma,G;X)
   = \begin{cases}
        2k^2 - 1    &(k: \text{even}),  \\
        2k^2 - 2    &(k: \text{odd}, \ \ge 3),  \\
        \quad\ \ 2  &(k=1).
     \end{cases}
$$
Thus, local rigidity (R) fails for any $k$ 
because the dimension of the deformation 
space is positive.
In the above formula,
 one sees that
 the dimension of the  deformation space 
$\mathcal{T}(\Gamma, G; X)$
has a different feature according to whether $k$ is even or odd. 
This will be explained 
by the criterion of properly discontinuous
 actions
 %our explicit description of $R(\Gamma, G; X)$ 
which involves the
existence of a non-zero real eigenvalue
 of a certain $k \times k$ matrix,
 whence the parity of $k$ counts.
Moreover, it follows from the complete
 description of 
$R(\Gamma, G; X)$  that
 we can determine for which 
$\varphi_0 \in R(\Gamma,G;X)$
 the stability (S) fails.

Our specific choice of $G$ was motivated by
 Lipsman's classification \cite{lipsman}
of maximal nilpotent affine transformation groups on
$\mathbb{R}^3$,
in which the two-step nilpotent group $G$ for $k=2$ played a crucial
 role. 
It is noteworthy that
 for any subgroup $\widetilde{G}$
 of the affine transformation group 
 $\operatorname{Aff}(\mathbb R^{k+1})$
 containing our specific $G$,
$R(\Gamma,\widetilde{G}; X)$ is not open in 
$\Hom(\Gamma,\widetilde{G})$ by Theorem~\ref{T:deformpara}, and consequently 
 both (R) and (S) fail.

\smallskip
The key idea of our proof is to take
 a connected subgroup $L$
   that contains $\Gamma$ as a cocompact
  discrete subgroup,
  and then to show that
  every injective
   homomorphism from $\Gamma$ into $G$
   extends uniquely to a continuous homomorphism
   from $L$ into $G$
   (an idea of {\it syndetic hull}).
A second step is to find explicitly
 $\Hom(L, G)$ in place of $\Hom(\Gamma, G)$,
 and to determine which homomorphism
 yields a properly discontinuous action.
Unlike the reductive case \cite{benoist, kobayashi89, kobayashi96},
 properly discontinuous actions for
  affine transformation groups
   on $\mathbb R^{k+1}$
   are still far from being fully understood in general,
  as one sees from the current status of the long-standing Auslander
   conjecture  (see \cite{xams} and
   references therein).
However, fortunately in our special setting,
 we can use the 
 criterion  \cite{nasrin}
  of proper actions for 
 two-step nilpotent Lie groups,
 which was obtained  as an 
 affirmative solution to
  Lipsman's conjecture \cite{lipsman}. 
Then, the final step to the proof of Theorem~\ref{T:deformpara}
(the description of $R(\Gamma, G; X)$)
is reduced to a certain problem of Lie algebras,
which we can solve explicitly.

\section{Description of deformation parameter}

This section gives a complete description
 of the parameter space
$R(\mathbb Z^{k}, G; \mathbb R^{k+1})$
of properly discontinuous $\mathbb{Z}^k$-actions on $\mathbb{R}^{k+1}$
 through a certain nilpotent affine
  transformation group $G$.
This is the first of the main results
 of this paper,
 and is stated in 
 Theorem~\ref{T:deformpara}.
Building on it,
we shall determine the deformation space 
$\mathcal{T}(\Gamma,G; \mathbb{R}^{k+1}) 
 \simeq R(\Gamma,G; \mathbb{R}^{k+1})/G$
in Section~\ref{sect:deformspace}
(see Theorem~\ref{T:deformspace}).
  
\subsection{Nilpotent affine transformation group} \label{SS:set1}

 We fix a positive integer $k$.
Our basic setting in this paper is:
\begin{align}\label{eqn:G}
\Gamma &:= \mathbb{Z}^k, 
\nonumber
\\ 
G &:= \left\{ \begin{pmatrix} 
I_k & \vec{x} & \vec{y}\\
0 & 1 & z \\
0 & 0 & 1
\end{pmatrix} :\, \vec{x}, \vec{y} \in \mathbb R^k, z \in \mathbb R \, 
\right\},
 \\
H &:= \left\{ \begin{pmatrix} 
I_k & \vec{x} & \vec{0} \\
0 & 1 & 0 \\
0 & 0 & 1
\end{pmatrix} :\, \vec{x} \in \mathbb {R}^k \, \right\}.
\nonumber
\end{align}
Then $G$ is a simply connected two-step nilpotent Lie group, and 
the homogeneous space $G/H$ is diffeomorphic to $\R^{k+1}$. 
We shall first give a description of $\Hom(\Gamma, G)$ 
in Proposition \ref{P:homomorphism}, and 
then
determine explicitly 
 $R(\Gamma, G; \mathbb{R}^{k+1})$
as its subset in Theorem \ref{T:deformpara}. 

Geometrically,
 this means that we determine
 all possible properly discontinuous affine actions of $\mathbb Z^k$
 on $\mathbb R^{k+1}$ preserving differential forms $d\xi_{k+1}$
and $d\xi_i \wedge d\xi_{k+1}$ $(1 \le i \le k)$,
where $(\xi_1, \ldots, \xi_{k+1})$ is the coordinate of
$\mathbb{R}^{k+1}$.

\subsection{Description of $\Hom (\Gamma, G)$} 
\label{SS:homomorphism} 

Any group homomorphism $\varphi : \Gamma \rightarrow G$ is determined by 
its
evaluation at generators of $\Gamma$. 
Taking a standard basis 
$\{e_1, \dots, e_k\}$ of the abelian group 
$\Gamma = \mathbb Z^k$, we regard $\Hom(\Gamma, G)$ as a subset of the 
direct product $\tG$ by the evaluation map: 
\begin{equation} \label{E:homo}
\Hom(\Gamma, G) \hookrightarrow \tG, \quad \! 
   \varphi \mapsto (\varphi(e_1), \dots, \varphi(e_k)).
\end{equation}

Let us describe the image of the injective map \eqref{E:homo}. 
For this, first we set
\begin{align}
M_1 &:= \{(\vec{x}, Y, \vec{z}) \in \mathbb R^k \oplus 
M(k, \mathbb R)\oplus \mathbb R^k : \vec{z} \neq \vec{0}\, \} 
\subset M(k, k+2; \mathbb R),  \label{D:a1} \\
M_2  &:=  M(k, 2k; \mathbb R).  \label{D:a2}
\end{align} 

\noindent
Then, $\dim M_1 = k(k+2)$ and $\dim M_2 = 2k^2$. 
Second, for $\ya{x}, \ya{y} \in \mathbb R^k$ and $z \in \mathbb R$,
 we define a $(k+2) \times (k+2)$ matrix by
\begin{equation} \label{E:set3}
  g(\ya{x}, \ya{y}, z) := \exp
\begin{pmatrix} 
     {\bigzero }_k & \ya{x} & \ya{y}\\
     0 & 0 & z \\
     0 & 0 & 0
\end{pmatrix}
=
\begin{pmatrix}
I_k & \ya{x} & \ya{y} + \frac{1}{2} z\ya{x}\\
0 & 1 & z \\
0 & 0 & 1
\end{pmatrix}.
\end{equation}

With this expression, our groups $G$ and 
$H$ (see Section~\ref{SS:set1}) are expressed as
\begin{align*}
   G &= \{ g(\ya{x}, \ya{y}, z) : 
   \ya{x}, \ya{y} \in \mathbb R^k, z \in \mathbb R \}, \\
   %\label{D:a4} \\
   H &= \{ g(\ya{x}, \ya{0}, 0): \ya{x} \in \mathbb R^k \}. 
   %\label{D:a5}
\end{align*}

Now we define maps 
$$\Psi_{i} \, :\, M_i \longrightarrow \tG \quad (i = 1, \, 2)$$ 
such that their $j$th ($1 \le j \le k$) components are respectively 
given by 
\begin{align} 
\label{E:e1}
\Psi_1(\ya{x}, Y, \ya{z})_j 
&:= g(z_j \ya{x}, \ya{y_j}, z_j) 
&&\text{for $Y= (\ya{y_1}, \cdots, \ya{y_k})$, $\trans{\ya z} = 
(z_1, \cdots, z_k)$,}
\\
\Psi_2(X, Y)_j 
&:= g(\ya{x_j}, \ya{y_j}, 0) 
&&\text{for $X = (\ya{x_1}, \cdots, \ya{x_k})$,
$Y= (\ya{y_1}, \cdots, \ya{y_k})$.}
\label{E:e2}
\end{align}

 %             0 & 0 & 0
  %   \end{pmatrix}
 %    \begin{pmatrix} I_k & z_j \ya{x} & \ya{y_j} + \frac12z_j^2 \ya{x}\\
  %            0 & 1 & z_j \\
   %           0 & 0 & 1
    % \end{pmatrix},
Then, the topological space
$\Hom(\Gamma,G)$ is described via \eqref{E:homo} as follows:

\begin{prop}[Description of $\Hom(\Gamma, G)$]
 \label{P:homomorphism}

\hfill\break
{\upshape 1)}\enspace
The maps $\Psi_1$ and $\Psi_2$ induce a bijection 
$$ 
\Psi_1 \cup \Psi_2 :
M_1 \cup M_2   % \underset{\Psi_1 \cup \Psi_2}
{\rarrowsim} 
\Hom(\Gamma, G).
$$ 
In particular, $\Psi_1$ and $\Psi_2$ are injective, and 
their images $\Psi_1 (M_1)$ and $\Psi_2 (M_2)$ are disjoint subsets 
contained in $\Hom(\Gamma, G)$.
\newline
{\upshape 2)}\enspace 
{\upshape(closure relation)}\enspace 
 $\Psi_2 (M_2)$ is closed in $\Hom(\Gamma, G)$, 
 whereas the closure of $\Psi_1(M_1)$ 
is given as
\begin{equation} \label{E:closure} 
\overline{\Psi_1 (M_1)} = \Psi_1 (M_1) \cup \Psi_2 (M^{d}_2).
\end{equation}
Here, $M^{d}_2$ is a subset of $M_2$ defined by 
\begin{equation} \label{E:m2d}
M^{d}_2 := \{ (X, Y) : X, \, Y \in M(k, \R), \; \rank X \le 1 \, \}.
\end{equation} 
\end{prop} 

We shall give a proof of this proposition
 in Section~\ref{S:3}.

\subsection{Description of $R(\Gamma, G; \mathbb{R}^{k+1})$} 
\label{SS:parameterspace}

Let us introduce the following subsets 
of $M_1$ and $M_2$: 
\begin{align} 
M^{r}_1 &:= \{(\ya{x}, Y, \ya{z}) \in M(k, k+2; \mathbb R) : 
\ya{z} \ne \ya{0}, \, \rank 
({}^t Y, {\ya{z}})
 =k \, \} 
\label{D:a3} 
\\
M^{r}_2 &:= \{(X, Y) \in M(k, 2k; \mathbb R) : 
\det(Y - \lambda X) \ne 0 \text{ for any $\lambda \in \R$} \, \}. 
 \label{D:a4}
\end{align} 

We are now ready to characterize $R(\Gamma, G; \mathbb{R}^{k+1})$ as a 
subset of $\Hom(\Gamma, G)$. 
Here is the first of the main results in
 this paper:

\begin{thm}[Description of $R(\Gamma,G;\mathbb{R}^{k+1})$]
 \label{T:deformpara} 
Let $G$ be a nilpotent Lie group defined as
\eqref{eqn:G} and
$\Gamma = \mathbb{Z}^k$.
Then,
the maps $\Psi_1$ and $\Psi_2$ (see \eqref{E:e1} and \eqref{E:e2}) 
induce the bijection 
$$
\Psi_1 \cup \Psi_2 :
 M^{r}_1 \cup M^{r}_2 \rarrowsim
R(\Gamma, G; \mathbb{R}^{k+1}).
$$
\end{thm}

We shall give a proof of this theorem in
Section~\ref{S:proof1}.

\subsection{Generic points of $R (\Gamma, G; \R^{k+1})$}

This subsection studies a generic part of $R (\Gamma, G; \R^{k+1})$
by analysing the sets $M_1^r$ and $M_2^r$ in detail.

\smallskip
For $X, \, Y \in M(k, \R)$, we define a polynomial of $\lambda$ by 
\begin{equation} \label{E:polynomial} 
f(X, Y ; \lambda) := \det (Y - \lambda X)      %% \label{D:a1}  
= \sum_{l = 0}^{k} a_l (X, Y) \lambda^{l}.
\end{equation} 
Here, we note 
\begin{align*}
a_0 (X, Y) &= \det Y, \\
a_{k-1} (X, Y) &= (-1)^{k-1} \sum_{i, j = 1}^{k} (-1)^{i+j} Y_{ij} 
\det \widehat{X}_{ij}, \\ 
a_{k} (X, Y) &= (-1)^k \det X, 
\end{align*}
where $\widehat{X}_{ij}$ denotes the submatrix obtained by deleting 
row $i$ and column $j$ from $X$. 

\smallskip
Let us recall from \eqref{D:a4} and \eqref{E:polynomial} that
$$M_2^{r} = \{(X, Y) \in M(k, 2k; \R) : f(X, Y ; \lambda) \ne 0 
\text{ for any } \lambda \in \R \, \}. $$ 
In order to give a `generic' part of
$R(\Gamma, G; \mathbb{R}^{k+1})$ by means of Theorem~\ref{T:deformpara},
 we set 
$$ M_{2}^{\operatorname{ro}} := 
\begin{cases} 
M_2^{r} \cap \{ (X, Y) : \det X \ne 0 \} & \text{($k$: even),} \\
M_2^{r} \cap \{ (X, Y) : \rank X = k-1, a_{k-1} (X, Y) \ne 0 \} 
& \text{($k$: odd).} 
\end{cases}
$$

\begin{prop}\label{P:2.4.1}\hfill\break
{\upshape1)}\enspace
$M_{1}^r$ is open dense in $M_1$. In particular, it has dimension $k(k+2)$. 

\noindent
{\upshape2)}\enspace
$ M_{2}^{\operatorname{ro}}$ is open in $M_{2}^r$, and the complement 
of $M_{2}^{\operatorname{ro}}$ in $M_{2}^r$ has a smaller 
dimension than that of $M_{2}^{\operatorname{ro}}$. The dimension of 
$ M_{2}^{\operatorname{ro}}$ is given by 
$$ \dim M_{2}^{\operatorname{ro}} = \begin{cases} 
2k^2 & (\text{$k$: even\/}), \\
2k^{2} - 1 & (\text{$k$: odd\/}). 
\end{cases}
$$ 
\end{prop}

\begin{proof}
1)\enspace  
Clear. 

2)\enspace  
For an even integer $l$, we consider a monic polynomial 
of the real variable $x$: 
$$g(x) = x^l + b_{l-1} x^{l-1} + \cdots + b_1 x + b_0. $$ 
Then, it attains its minimum,
 denoted by $m ( b_0, b_1, \cdots, b_{l-1})$, 
which is a continuous function of the real coefficients 
$b_0, b_1, \cdots, b_{l-1}$. 

\medskip
\noindent
{\bf Case 1 \, ($k$ : even)} \, 
Suppose  
$(X, Y) \in M_{2}^{\ro}$. 
Then, the monic polynomial 
$ \frac{f(X, Y; \lambda)}{a_k (X, Y)}$ must be positive for all 
$ \lambda \in \R$. 
Thus we have
$$M_{2}^{\ro} = \left\{ (X, Y) \in M_2 : \det X \ne 0, \,
m \left( \frac{a_0 (X, Y)}{a_k (X, Y)}, \cdots, 
\frac{a_{k-1} (X, Y)}{a_k (X, Y)} \right) > 0 \right\}.$$ 
Hence, $M_{2}^{\ro}$ is open in $M_2 = M(k, 2k; \R)$. 

\smallskip
To see $M_{2}^{\ro} \ne \emptyset$,
 we set $J_k := \begin{pmatrix} 
\bigzero & \raisebox{.7ex}{$-I_{\frac{k}{2}}$} \\ 
I_{\frac{k}{2}} & \bigzero 
\end{pmatrix} \in M(k, \R).$
 Then $(J_k, \, I_k) \in M_{2}^{\ro}$
  because 
$$ f (J_k, I_k; \lambda) = \left(\det \begin{pmatrix} 
1 & \lambda \\ 
- \lambda & 1 
\end{pmatrix} \right)^{\frac{k}{2}} = 
(1+ \lambda^{2})^{\frac{k}{2}} 
> 0.$$ 
Hence, $M_{2}^{\ro} \ne \emptyset$ and $\dim M_{2}^{\ro} = 
\dim M_2 = 2k^2.$

\smallskip
Next, suppose $(X, Y) \in M_2^{r} \setminus M_{2}^{\ro}$. Then 
$\det X = 0$ by definition. Furthermore, it follows from 
$f(X, Y; \lambda) \ne 0$ 
for any $\lambda \in \R$ that the coefficient $a_{k-1} (X, Y)$ 
of $\lambda^{k-1}$ must vanish because 
$f(X, Y; \lambda) = a_{k-1} (X, Y) \lambda^{k-1} + \cdots + 
a_{0} (X, Y)$ and $k-1$ is odd. 
Thus we have seen that 
$$M_2^{r} \setminus M_{2}^{\ro} \subset \{ (X, Y) \in M_2 : 
\det X = a_{k-1} (X, Y) = 0 \, \}.$$ 
Hence the complement $M_2^{r} \setminus M_{2}^{\ro}$ has at least codimension 
two in $M_2^{r}$. 

\medskip
\noindent
{\bf Case 2 ($k$: odd)} First, we claim 
$$M_2^{r} \subset \{ (X, Y) \in M(k, 2k; \R) : \det X = 0 \, \}.$$ 
In fact, since $k$ is odd,
the polynomial $f(X, Y; \lambda)$
 of the real variable $\lambda $ has zeros 
unless the top term 
$a_{k} (X, Y) \lambda^{k}$ vanishes.
Therefore, $a_{k} (X, Y) \; (= (-1)^k \det X) = 0$ if 
$(X, Y) \in M_2^{r}$. 

\smallskip
Next, let us prove that $M_{2}^{\ro}$ is open in the set 
\begin{align*}
S &:= \{ (X, Y) \in M(k, 2k; \R) : \det X =0, \, \operatorname{grad} 
\det X \ne 0 \, \} \\ 
&= \{ (X, Y) \in M(k, 2k; \R) : \rank X = k-1 \, \}.
\end{align*} 

Suppose $(X, Y) \in M(k, 2k; \R)$ satisfies 
$$\det X =0 \text{ and } a_{k-1} (X, Y) \ne 0.$$ 
Then $f(X, Y; \lambda) \ne 0$ for any $\lambda \in \R$ if and only if 
the monic polynomial 
$$\frac{f(X, Y; \lambda)}{a_{k-1} (X, Y)} = \lambda^{k-1} + 
\sum_{i=0}^{k-2} \frac{a_{i} (X, Y)}{a_{k-1} (X, Y)} \lambda^i $$ 
is positive for all $\lambda$. Thus, we have seen 
\begin{align*} 
M_{2}^{\ro} &= \left\{ (X, Y) \in M(k, 2k; \R) : \begin{array}{cl} 
& \rank X = k-1, \, a_{k-1} (X, Y) \ne 0, \\ 
& m \left( \frac{a_0 (X, Y)}{a_{k-1} (X, Y)}, \cdots,  
\frac{a_{k-2} (X, Y)}{a_{k-1} (X, Y)} \right) > 0. 
\end{array} \right\}. 
\end{align*}

\smallskip
It is now clear that $M_{2}^{\ro}$ is open in $S$. 
To see $M_{2}^{\ro} \ne \emptyset$,
we set $J'_k:= \begin{pmatrix} 
J_{k-1} & \bigzero  \\ 
&& \\
\bigzero & 0 \\ 
\end{pmatrix} \in M(k; \mathbb R)$.
Then,  
$( J'_k, I_k) \in M_{2}^{\ro}$
because 
$$ f \left( J'_k, I_k; \lambda \right) = (1+ \lambda^{2})^{\frac{k-1}{2}} 
> 0.$$ 
Hence, $M_{2}^{\ro}$ is a non-empty open subset of $S$. 
 Since $S$ is a non-singular 
manifold of dimension $2k^2 - 1$, so is $M_{2}^{\ro}$.

\smallskip
Finally, it follows from the definition of $M_{2}^{\ro}$ that 
$M_2^{r} \setminus M_{2}^{\ro}$ is contained in the algebraic variety: 
$$\{(X, Y) : \det X = 0, \, a_{k-1} (X, Y)
\|\operatorname{grad} \det X \|^2 = 0 \, \},$$ 
which is  of dimension $2 k^2 -2$. Thus, 
Proposition~\ref{P:2.4.1} has been proved. 
\end{proof}

\begin{rem}
Proposition~\ref{P:2.4.1} implies that 
$ M^{r}_2$ contains an open subset of $M_2$ if and only if $k$ is even. 
However, $M^{r}_2$ itself is not open in $M_2$ 
even if $k$ is even because $(O, I_k) \in M^{r}_2$
 is not an inner point,
 as we shall see in the following example:
\end{rem}

\begin{exa} \label{R:halfdim}
Take a half dimensional affine subspace
$$
  V = \{(X, I_k): X \in M(k, \mathbb R) \}
$$
 of $M(k, 2k; \mathbb R).
$
Then we have
$$
V \cap M^{r}_2 = \{(X, I_k): \text{any eigenvalue of $X$ 
is in $\mathbb C\setminus \mathbb R^\times$} \}.
$$
This gives a partial information on $R(\Gamma, G; \mathbb{R}^{k+1})$, 
and was proved in Lipsman \cite[Theorem 4.4]{lipsman} 
 for $k= 2$ as a crucial step to the classification of maximal
nilpotent affine subgroups that act properly on $\mathbb{R}^3$,
and was generalized
 in \cite{nasrin1} for $k\ge3$. Our proof here
 is different and simpler. 
\end{exa}

\section{Description of $\Hom(\Gamma, G)$}
\label{S:3}

This section determines
 $\Hom(\Gamma, G)$ explicitly,
  and gives a proof of Proposition~\ref{P:homomorphism}.

\subsection{Parametrization of $\Hom(\Gamma,G)$} \label{SS:parameter} 

Recall from \eqref{E:homo} that any $\varphi \in \Hom(\Gamma, G)$ 
is determined by $\varphi(e_j)$  $(1 \le j \le k)$, which we write 
as $\varphi(e_j) = g(\ya{x_j}, \ya{y_j}, z_j)$ for some 
$\ya{x_j}, \, \ya{y_j} \in \R^k$ and $z_j \in \R$ according to 
\eqref{E:set3}. Collecting these data
$\ya{x_j}, \ya{y_j}$ $(1 \le j \le k)$ and 
$\ya{z} := \trans(z_1, \dots, z_k)$,
 we obtain an injective map defined by
\begin{equation} \label{E:homo1}
\Hom(\Gamma, G) \to M(k, 2k+1; \mathbb R), \, \ 
   \varphi \mapsto %\begin{pmatrix}
                    \bigl(    \ya{x_1}, \dots, \ya{x_k}
                   ;     \ya{y_1}, \dots, \ya{y_k}
                   ; \ya{z} \bigr).
\end{equation}

\smallskip 
Let us determine the image of \eqref{E:homo1}. 
Since $\varphi \in \Hom(\Gamma, G)$ satisfies
 $$\varphi(e_i) \varphi(e_j) = \varphi(e_i + e_j) = \varphi(e_j) \varphi(e_i)$$
 for any $i, j$  $(1 \le i, j \le k)$, we have
\begin{equation} \label{E:equn}  
        z_i \ya{x_j} = z_j \ya{x_i}
   \quad
     \text{ for any $i, j$ $(1 \le i, j \le k)$}.
\end{equation}
Conversely, given $k$ elements $g_1, \cdots, g_k$ in $G$ that mutually commute, we can 
define a group homomorphism 
$\varphi: \Gamma \to G$ by $\varphi ( \sum_{j= 1}^{k} m_j e_j ) := 
g_{1}^{m_1} \cdots g_{k}^{m_k}$. 
Therefore, the image of (\ref{E:homo1}) is 
characterized by the condition \eqref{E:equn}, that is, 
we have a bijection:
$$
\Hom(\Gamma,G) \simeq
   \{
              \begin{pmatrix}
                        \ya{x_1}, \dots, \ya{x_k};
                    Y; \ya{z} %                  \\      z_1 & \dots & z_k
 \end{pmatrix}: Y \in M(k, \mathbb R);
 \ya{x_1}, \dots, \ya{x_k}, \ya{z} \ 
 \text{satisfies (\ref{E:equn})} \}.
$$
We shall find all solutions of (\ref{E:equn}), 
according to the following two cases:
(a)\; $\ya{z} \neq \ya{0}$ and (b)\; $\ya{z} = \ya{0}$.

\smallskip
In the case (a),
           there exists uniquely an element $\ya{r} \in \mathbb R^k$
         such that $\ya{x_j} = z_j \ya{r}$ for all $j$ ($1 \le j \le k$).
This amounts to  $\Psi_1(M_1)$ (see \eqref{E:e1} for the definition 
of $\Psi_1$).

In the case (b),
           any $\ya{x_1}, \dots, \ya{x_k}$ solves \eqref{E:equn}.
This amounts to $\Psi_2(M_2)$ (see \eqref{E:e2} for the 
definition of $\Psi_2$).

\smallskip
Hence, $\Hom (\Gamma, G)$ is the disjoint union of $\Psi_1(M_1)$ 
and $\Psi_2(M_2)$. Thus we have completed the proof of 
Proposition~\ref{P:homomorphism}~(1). 

\qed 

\subsection{Closure relation in $\Hom(\Gamma,G)$}

This subsection gives a proof of Proposition~\ref{P:homomorphism}~(2).
It is clear that $\Psi_2(M_2)$ is a closed set.
Let us consider the closure of $\Psi_1(M_1)$,
 and find its boundary.
 What we need is to prove: 
$$
\overline{\Psi_1 (M_1)} \cap \Psi_2 (M_2) = \Psi_2 (M^{d}_2).
$$
{\bf Proof of the inclusion $\supset$ :} \, 
Suppose $(X, Y) \in \Psi_2 (M^{d}_2)$. Since $\rank X \le 1$, we find 
$\ya{x} \in \R^k$ and $\ya{a} \in \R^k \setminus \{0\}$ such that 
$X = \ya{x} \trans{\ya{a}}.$ 
In light of the obvious formula
$$g (a_j \ya{x}, \ya{y_j}, 0) = \lim_{l \to \infty} 
g \bigl( \frac{a_j}{l} l \ya{x}, \ya{y_j}, \frac{a_j}{l} \bigr)  
\quad \! (1 \le j \le k),$$ 
we conclude from the definitions \eqref{E:e1} and \eqref{E:e2} 
of $\Psi_1$ and $\Psi_2$ that 
$$\Psi_2 (X, Y) = \lim_{l \to \infty} \Psi_1 \bigl( l \ya{x}, Y, 
\frac{{\ya{a}}}{l} \bigr).$$

As $(l \ya{x}, Y, \frac{{\ya{a}}}{l})$  $(l = 1, 2, \cdots)$ is a 
sequence of $M_1$, we have proved the inclusion $\supset$. 
\newline 
{\bf Proof of the inclusion $\subset$ :} \, 
Take any sequence $(\ya{x^{(l)}}, Y^{(l)}, 
\ya{z^{(l)}})$ in $M_1$ such that 
$\Psi_1 (\ya{x^{(l)}}, Y^{(l)}, \ya{z^{(l)}})$ converges to an element of 
$\Psi_2 (M_2)$, say, $\Psi_2 (X, Y)$ for some $X, \, Y \in M(k, \R)$. 
Then the formula 
$$\lim_{l \to \infty} \Psi_1 (\ya{x^{(l)}}, Y^{(l)}, \ya{z^{(l)}}) = 
\Psi_2 (X, Y)$$ 
implies that $X$ is the limit of 
$X^{(l)}:= 
(z_1^{(l)} \ya{x^{(l)}}, \cdots, z_k^{(l)} \ya{x^{(l)}})$ as $l$ 
tends to 
infinity. Since $\rank{X^{(l)}} \le 1$, its limit also satisfies 
$\rank X \le 1$. Thus we have proved the inclusion $\subset$. 

Thus, Proposition~\ref{P:homomorphism}~(2)
 is proved.
\qed

\section{Proof of Theorem~\ref{T:deformpara}} \label{S:proof1}
This section gives a proof of
Theorem~\ref{T:deformpara}.
Our strategy here is to rewrite the condition of
$R(\Gamma,G; \mathbb{R}^{k+1})$, in particular,
the condition for properly discontinuous actions in the following scheme:
$$
\begin{array}{llll}
&  \Gamma  &\text{discrete subgroup}
&\text{(see \eqref{eqn:R})}
\\
\Rightarrow& L = \overline{\Gamma} &\text{its syndetic hull}
&\text{(see Proposition~\ref{L:coro})}
\\
\Rightarrow& \mathfrak{l} &\text{its Lie algebra.}
&\text{(see Section~\ref{SS:reform})}
\end{array}
$$

\subsection{Proper actions and properly discontinuous actions} 
\label{SS:proper} 

In dealing with properly discontinuous actions of a discrete group,
a more general notion ``proper action'' is sometimes useful.
We recall:
\begin{defn}[Palais \cite{palais}]
Suppose that a locally compact topological group $L$ 
acts continuously on a Hausdorff,
locally compact  space $X$. For a subset $S$ 
of $X$, we define a subset of $L$ by 
$L_{S} = \{ \gamma \in L : \gamma S \cap S \ne \emptyset \}.$
The $L$-action on $X$ is said to be \textit{proper}
 if $L_{S}$ is compact
 for every compact subset $S$ of $X$. 
\end{defn} 

\smallskip
We note that the $L$-action is \textit{properly discontinuous} if 
$L$ is a discrete group and if the $L$-action is proper.

\smallskip
The following elementary observation is a bridge between the action of a 
discrete group and that of a connected group.

\begin{obs}[{\cite[Lemma 2.3]{kobayashi89}}]
 \label{O:bridge}
Suppose a locally compact group $L$ acts on a Hausdorff, locally compact 
space $X$. Let $ \Gamma$ be a cocompact discrete subgroup of $L$. Then
\begin{enumerate}
\renewcommand{\labelenumi}{\upshape\theenumi)}
   \item The $L$-action on $X$ is proper if and only if the $ \Gamma $-action 
   is properly discontinuous.
   \item $L \backslash X$ is compact if and only if\/
 $ \Gamma \backslash X$ is 
   compact.
\end{enumerate}
\end{obs}

\smallskip 
\subsection{Extension from a discrete subgroup} 

Suppose we are in the setting of Section~\ref{SS:set1}. 
We set 
$$
       L := \mathbb R^k
$$
 and regard $\Gamma =\mathbb Z^k$ as a cocompact discrete subgroup of $L$. 
  %\smallskip
We write $\Hom(L, G)$ for the set of 
continuous group homomorphisms from $L$ into $G$. 
In our setting \eqref{eqn:G},
every homomorphism from $\Gamma$ into $G$ extends uniquely to a
  continuous homomorphism from $L$ to $G$.
That is, we have:
\begin{lem} \label{L:rest}
The restriction map\/ 
$\Hom(L, G) \to \Hom(\Gamma, G), \psi \mapsto \psi|_\Gamma$
is bijective.
\end{lem}

\begin{proof}
As $G$ is a simply connected nilpotent group, the 
exponential map, $\exp : \mathfrak g \rightarrow G$ is bijective. 
We write $\log$ for its inverse. Then,
$\psi \in \Hom(L,G)$ satisfies
\begin{equation}\label{eqn:res}
\psi(\sum_{j=1}^k a_j e_j) = \exp \bigl(\sum_{j=1}^k a_j \log 
\psi(e_j)\bigr) \text{ for any } a_1, \cdots, a_k \in \R.
\end{equation}
This shows that the homomorphism $\psi$ is determined by 
its restriction
$\psi|_\Gamma$.
Conversely, the formula \eqref{eqn:res} also indicates how to extend a
homomorphism from $\Gamma$ to $L$.
Thus we have proved Lemma \ref{L:rest}.
\end{proof}

In light of Lemma \ref{L:rest},
any property of $\psi$ should be expressed in terms of the restriction
$\psi|_{\Gamma}$ in principle.
We show:

\begin{lem} \label{L:disc} 
The following two conditions on
 $\psi \in \Hom(L, G)$ are equivalent:
\begin{newenumerate}
\item $\psi$ is injective.
\item $\psi|_{\Gamma}$ is injective and $\psi(\Gamma)$ is discrete in $G$.
\end{newenumerate}
\end{lem}

\begin{proof}
Since $G$ is a simply connected nilpotent Lie group, any connected 
subgroup of $G$ is closed. Therefore, $\psi : \, L/\Ker \psi \to G$ is a 
homeomorphism onto a closed subgroup of $G$.
In particular, $\psi(\Gamma)$ is discrete in $G$
 if and only if $\Gamma/\Gamma \cap \Ker \psi$ is discrete in $L/\Ker \psi$.
Now, it is clear that (i) implies (ii). 

Conversely, if $\psi$ is not injective and if $\psi|_{\Gamma}$ is injective,
 then the composition map $\Gamma \subset L \to L/\Ker \psi$
 is injective with non-discrete image because 
 $\rank \Gamma < \dim (L/\Ker \psi)$. 
Hence, $\psi(\Gamma)$ is not discrete in $G$, too.
Thus, (ii) also implies (i).
\end{proof}

\subsection{A continuous analogue of properly discontinuous actions} 

Following Observation \ref{O:bridge},
 we amplify Lemma \ref{L:disc} 
with the condition of proper actions on the homogeneous space $G/H$:

\begin{lem} \label{L:equiv}
Let $\psi \in \Hom(L, G)$ and $\varphi = \psi|_{\Gamma}$ 
(see Lemma \ref{L:rest}).
Then the following two conditions are equivalent:
\begin{newenumerate}
\item $\psi: L \to G$ is injective and $\psi(L)$ acts properly on $G/H$.
\item $\varphi: \Gamma \to G$ is injective and $\varphi(\Gamma)$ 
 acts properly discontinuously and freely on $G/H$.
\end{newenumerate}
\end{lem}

\begin{proof}
(i) $\Rightarrow$ (ii):
Since $\psi$ is injective, it follows from Lemma \ref{L:disc} that 
$\varphi(\Gamma)$ is discrete in a closed subgroup $\psi(L)$.
Therefore, $\varphi(\Gamma)$ acts properly discontinuously on $G/H$ 
because $\psi(L)$ acts properly on $G/H$.
Furthermore, any properly discontinuous action of $\varphi(\Gamma)$ is 
automatically free because $\varphi(\Gamma) \simeq \Gamma$ is torsion-free.
Hence (ii) is proved.

\smallskip
(ii) $\Rightarrow$ (i):
If $\varphi(\Gamma)$ acts properly discontinuously on $G/H$
 then $\varphi(\Gamma)$ is discrete in $G$.
Hence, $\psi : L \to G$ is injective by Lemma \ref{L:disc}.
Furthermore, $\psi(L)$ with its relative topology contains $\varphi(\Gamma)$ as a cocompact discrete subgroup.
Therefore, $\psi(L)$ acts properly on $G/H$ by Observation 
\ref{O:bridge} (1).
Thus, we have proved
 the implication (ii) $\Rightarrow$ (i).
\end{proof}

We are ready to 
characterize $R(\Gamma,G; \mathbb{R}^{k+1})$ 
 by means of the connected
subgroup $L$:

\begin{prop} \label{L:coro}
Under the isomorphism
$\Hom(L,G) \rarrowsim \Hom(\Gamma,G)$ 
in Lemma~\ref{L:rest},
we have
\begin{align*}
R(\Gamma, G; \mathbb{R}^{k+1}) \simeq \{\psi \in \Hom(L, G): \ 
&\text{{\upshape i)}\ 
$\psi$ 
 is injective},
\\
&\text{{\upshape ii)}\ 
 $\psi(L)$ acts properly on $G/H$} \}.
\end{align*}
\end{prop}

\subsection{Reformulation of $R(\Gamma, G; \mathbb{R}^{k+1})$} 
\label{SS:reform}

So far, we have transferred proper discontinuity and freeness of
discrete group actions into a certain property of connected group
actions. 
Now, let us rewrite the latter condition  in terms of 
Lie algebras. We use the German lower case letters $ \mathfrak g$, $ \mathfrak h$ 
and $ \mathfrak l$ to denote the Lie algebras of $G$, $H$ and $L$ 
respectively. 
We write $d \psi$ for the differential of $\psi \in \Hom(L, G)$.
Consider the following conditions on $d\psi$:
\begin{gather}
d \psi : \mathfrak l \to \mathfrak g \text{ is injective}, \label{E:diff1} \\
d \psi(\mathfrak l) \cap \bigcup_{g \in G} \Ad(g) \mathfrak h = \{0\}. \label{E:diff2}
\end{gather}

\smallskip
Now we can restate Proposition~\ref{L:coro} as 
\begin{prop} \label{P:isom}
Under the isomorphism $\Hom(L, G) \rarrowsim \Hom(\Gamma, G)$ 
(see Lemma \ref{L:rest}), we have
$$
R(\Gamma, G; \mathbb{R}^{k+1}) \simeq \{ \psi \in \Hom(L, G) \, : d\psi 
\text{ satisfies \eqref{E:diff1} and \eqref{E:diff2}} \, \}.
$$
\end{prop}

\begin{proof} 
Any connected subgroup of a simply connected nilpotent Lie group is 
simply connected. Hence, $d \psi$ is injective if and only if $\psi$ 
is injective. Now use the criterion of proper actions for a 
homogeneous space of a two-step nilpotent Lie group $G$ 
as follows. 
\end{proof} 

\begin{lem}[{\cite[Theorem 2.11]{nasrin}}] 
\label{L:msthesis}
Let $G$ be a simply connected Lie group,
and $H, L$ its closed subgroups.
Suppose $G$ is a two-step nilpotent Lie group,
which means that the commutator subgroup of $G$ is contained in the
centre of $G$.
Then, the
 following three conditions on $\psi$ are equivalent.
\begin{newenumerate}
\item $\psi(L)$ acts on $G/H$ properly.
\item $\psi(L) \cap gHg^{-1} \, = \, \{e\} \; \text{for all} \; g \in G.$
\item $d \psi(\mathfrak l) \cap \bigcup_{g \in G} \Ad(g) \mathfrak h
 = \{0\}.$
\end{newenumerate}
\end{lem}

\begin{rem}
Lemma~\ref{L:msthesis} gives an affirmative solution to Lipsman's
conjecture \cite{lipsman} for two-step nilpotent Lie groups.
Recently,
Baklouti--Khlif \cite{baklouti}
and Yoshino \cite{yoshino3step}
proved independently that Lipsman's conjecture is still true for
three-step nilpotent Lie groups.
\end{rem}

\subsection{Completion of the proof of Theorem \ref{T:deformpara} }
\label{SS:pf2}

Now let us complete the proof of Theorem \ref{T:deformpara}. 
We have already reduced it to a problem of Lie algebras.
Now, we use the following: 

\begin{lem} \label{L:va}
We define the variety in $\mathfrak{g}$ by 
$$\mathcal{V} = \bigcup_{g \in G} \Ad(g) \mathfrak h.$$ 
Then we have 
\begin{eqnarray} \label{E:variety}
\mathcal{V} &=& \{ W - [W,V] : W \in \mathfrak{h}, V \in \mathfrak{g} \} \notag \\
         &=& \{ \begin{pmatrix}
         \bigzero & \ya{x} & b \ya{x} \\
         0 & 0 & 0 \\
         0 & 0 & 0
         \end{pmatrix} 
          \; : \ya{x} \in \mathbb R^k, b 
         \in \mathbb R \}. 
\end{eqnarray}
\end{lem}

\begin{proof}
Elementary computation.
\end{proof} 

According to the parametrization
$$
\Psi_1 \cup \Psi_2 :
M_1 \cup M_2 
\stackrel{\sim}{\to}
\Hom(L,G) \stackrel{\sim}{\to} \Hom(\Gamma, G)
$$
given in Proposition~\ref{P:homomorphism} and Lemma~\ref{L:rest}, we
examine if $d\psi$ satisfies \eqref{E:diff1} and \eqref{E:diff2} for
$\psi \in \Psi_1(M_1)$
and
$\psi \in \Psi_2(M_2)$, 
respectively.
The following proposition gives criteria for 
\eqref{E:diff1} and \eqref{E:diff2}:

\smallskip
\begin{prop}\label{P:comp}\hfill\break
{\upshape1)}\enspace
Let $\psi := \Psi_1(\ya{x}, Y, \ya{z})$ 
for $(\ya{x}, Y, \ya{z}) \in M_1$. 
Then,  we have the following equivalence:
$$
\psi 
\text{ satisfies \eqref{E:diff1} } \Longleftrightarrow 
\rank
({}^t Y, {\ya z} )
 = k. 
 $$ 
In this case,
 $\psi$ satisfies \eqref{E:diff2}, too. 

\noindent
{\upshape2)}\enspace
Let $\psi := \Psi_2(X, Y)$ for $(X, Y) \in M_2$. 
Then 
$$\psi 
\text{ satisfies \eqref{E:diff1} } \Longleftrightarrow 
\rank \begin{pmatrix} 
X \\
Y 
\end{pmatrix} = k. $$ 
In this case,
we have the following equivalence:
$$
\psi 
\text{ satisfies \eqref{E:diff2} } \Longleftrightarrow 
\det (Y - b X) \ne 0 \text{ for any $b \in \R $.}$$ 
\end{prop}

\smallskip
\begin{proof}
1)\enspace
 It follows from the definition of $\Psi_1$ 
(see \eqref{E:e1}) that 
\begin{equation} \label{E:dpsi1}
 d \psi (\ya{a}) = \begin{pmatrix}
         \bigzero & \langle \ya{a}, \ya{z} \rangle \ya{x} & Y \ya{a} \\
         0 & 0 & \langle \ya{a}, \ya{z} \rangle \\
         0 & 0 & 0
         \end{pmatrix}.
\end{equation}
Here, $\langle \ , \ \rangle$ denotes the standard inner product on
$\mathbb{R}^k$.
Then 
\begin{align*}
\begin{split}
\psi \text{ satisfies } \eqref{E:diff1} &\Longleftrightarrow \; \{ \ya{a} 
\in \mathfrak{l} \, : \, \langle \ya{z}, \ya{a} \rangle = 0, 
\, Y \ya{a} = \vec{0} \, \} = \{0\} \\
&\Longleftrightarrow \; \rank 
({}^t Y, {\ya{z}} )
 = k. 
\end{split}
\end{align*} 

Furthermore, by \eqref{E:variety} and \eqref{E:dpsi1} 
$$
d \psi(\ya{a}) \in \mathcal{V} \Longleftrightarrow 
d \psi(\ya{a}) = \vec{0}.
$$ 
Therefore, \eqref{E:diff1} implies \eqref{E:diff2}.

\smallskip
2)\enspace  
It follows from the definition of $\Psi_2$ (see \eqref{E:e2}) that
\begin{equation} \label{E:dpsi2}
 d \psi (\ya{a}) = \begin{pmatrix}
         \bigzero & X \ya{a}  & Y \ya{a} \\
         0 & 0 & 0 \\
         0 & 0 & 0
         \end{pmatrix}.
\end{equation} 
Then, 
\begin{align*}
\begin{split}
\psi \text{ satisfies } \eqref{E:diff1} &\Longleftrightarrow  \; 
\{ \ya{a} \in \mathfrak{l} \, : 
\, X \ya{a} = \vec{0}, \, Y \ya{a} = \vec{0} \, \} = \{\vec{0}\} \\
&\Longleftrightarrow \; \rank 
\begin{pmatrix} 
X \\
Y 
\end{pmatrix} = k.
\end{split}
\end{align*}

\smallskip
Suppose \eqref{E:diff1} is satisfied. It follows from 
\eqref{E:variety} and \eqref{E:dpsi2} that 
\begin{align*}
&\psi \text{ satisfies } \eqref{E:diff2} 
\\
\Longleftrightarrow \; 
&\text{there is no 
$b \in \R$ and $\ya{a} \in \R^k$ such that $Y \ya{a} =b X \ya{a} \ne \vec{0}$} \\ 
\Longleftrightarrow \; 
&\det(Y- b X) = 0 
 \text{ has no real solution for any $b \in \R$.}
\end{align*}
Hence, Proposition is proved.
\end{proof}

\smallskip
We note that if
$\det(Y-bX) \ne 0$
for any $b \in \mathbb{R}$ then
$\rank \begin{pmatrix} X\\Y \end{pmatrix} = k$
because $\det Y \ne 0$.
Then, Theorem \ref{T:deformpara} follows from Propositions \ref{P:isom} 
and \ref{P:comp}. 
Hence, we have completed the proof of Theorem \ref{T:deformpara}. 
\qed

\section{Deformation space} 
\label{sect:deformspace}

Building on the description of
 the parameter space $R(\Gamma, G; \mathbb R^{k+1})$
  of discontinuous groups given in
  Theorem~\ref{T:deformpara},
   we determine explicitly the deformation space
  $\mathcal{T} (\Gamma, G; \mathbb{R}^{k+1})$.
This is stated in Theorem~\ref{T:deformspace} and
Corollary~\ref{cor:dimdefo}, 
and is
the second of the main results
   of this paper. 
\subsection{Description of the deformation space 
  $\mathcal{T} (\Gamma, G; \mathbb{R}^{k+1})$} 

We define subsets of $M^{r}_{1}$ and $M^{r}_{2}$, respectively, by 

\begin{align} \label{D:a5}
D_1 &:= \left\{ \bigl( \ya{x}, 
{}^t (
{\ya{\eta_1}} 
,\ldots, 
{\ya{\eta_k}} 
), 
 \ya{z} \bigr) 
\in M(k,k+2;\mathbb{R})
: 
\begin{array}{cl}
1) & \ya{z}\neq \ya{0}, \\
2) & \vec{\eta_j} \perp \vec{z} \; \; 
(1 \le j \le k), \\
3) & \rank (\vec{\eta_1}, \cdots, \vec{\eta_k}) = k-1.
\end{array} 
\right\},
\end{align}

\begin{align} \label{D:a6}
D_2
 & := \left\{(X, Y) \in M(k, \mathbb R)\oplus M(k, \mathbb R): 
 \begin{array}{cl} 
1) & \operatorname{Trace} (X \trans{Y}) = 0,  \\
2) & \det(Y - \lambda X) \ne 0 \text{ for any $\lambda \in \R$.} 
\end{array} 
\right\}. 
\end{align} 

We note that the third condition in \eqref{D:a5} asserts that 
$\rank (\vec{\eta_1}, \cdots, \vec{\eta_k})$
 attains its maximum 
because all the vectors $\vec{\eta_j}$ $(1 \le j \le k)$
are orthogonal to $\vec{z}$.

\medskip

We retain the setting as in Section~\ref{SS:set1}.
In particular,
$\Gamma = \mathbb{Z}^k$ and $G$ is a nilpotent affine transformation
group defined in \eqref{eqn:G}.
For $i = 1, \, 2$, we denote by 
$$
\overline{\Psi}_{i} : M_{i}^{r} \rightarrow 
\mathcal{T} (\Gamma, G; \mathbb{R}^{k+1})
$$ 
the composition of 
$\Psi_{i} : M_{i}^{r} \rightarrow R(\Gamma, G; \mathbb{R}^{k+1})$ 
(see \eqref{E:e1} and \eqref{E:e2})
and 
the natural quotient map 
$R(\Gamma, G; \mathbb{R}^{k+1}) \rightarrow 
\mathcal{T} (\Gamma, G; \mathbb{R}^{k+1}).$ 

\medskip

Here is an explicit description
 of the deformation space:

\begin{thm} \label{T:deformspace} 
The maps $\overline{\Psi}_{1}$ and $\overline{\Psi}_{2}$ induce the 
following bijection:
$$
\overline{\Psi}_{1} \cup \overline{\Psi}_{2}:
D_1 \cup D_2   {\rarrowsim} 
\mathcal{T} (\Gamma, G; \mathbb{R}^{k+1}).$$
\end{thm}

In particular,
we find the dimension of the deformation space:

\begin{corollary}\label{cor:dimdefo}
The deformation space $\mathcal{T}(\Gamma,G;\mathbb{R}^{k+1})$ 
contains a smooth manifold $\mathcal{T}'$ as its open dense subset,
where the dimension of $\mathcal{T}'$ is given by
$$
\dim \mathcal{T}' =
   \begin{cases}
      2k^2 - 1 & (k: \text{even}),
   \\
      2k^2 - 2 & (k: \text{odd}, \ \ge 3),
   \\
      2        & (k = 1).
   \end{cases}
$$
\end{corollary}

\subsection{Proof of Theorem \ref{T:deformspace} and
   Corollary~\ref{cor:dimdefo} }
We let $G$ act on $M(k, k+2; \R)$ and $M(k, 2k; \R)$, 
respectively, as follows: 
\newline 
for $h = \begin{pmatrix} 
I_k & \vec{a} & \vec{b}\\
0 & 1 & c \\
0 & 0 & 1
\end{pmatrix} \in G$,
the actions of $h$ are given by
\begin{align} \label{E:e3}
(\ya{x}, Y, \ya{z}) 
&\mapsto (\ya{x}; 
Y + (\vec{a} - c \vec{x})\trans{\ya{z}}, \ya{z}) 
&&\text{on\quad $M(k, k+2; \R)$},
\\
\label{E:e4}
(X, Y) 
&\mapsto (X, Y - cX)% \; 
&&\text{on\quad $M(k, 2k; \R)$}. 
\end{align}

\begin{prop}\label{P:keyprop}\hfill\break
{\upshape1)}\enspace 
Both $M_1$ and $M^{r}_1$  are $G$-stable subsets of $M(k, k+2; \R)$. 

\noindent
{\upshape2)}\enspace 
$M^{r}_2$ is a $G$-stable subset of $M_2 = M(k, 2k; \R).$ 

\noindent
{\upshape3)}\enspace 
For $i = 1, \, 2$, the maps 
$\Psi_{i} : M^{r}_i \rightarrow \Hom(\Gamma, G)$ respect $G$-actions. 

\noindent
{\upshape4)}\enspace 
For $i = 1, \, 2$, $D_i$ are complete representatives of the 
$G$-orbit on $ M^{r}_i$.
\end{prop}

\begin{proof}
1) Clear from the definitions \eqref{D:a1} and \eqref{D:a3} 
of $M_1$ and $M^{r}_1$. 
\newline 
2) Clear from the definition \eqref{D:a4} of $M^{r}_2$. 
\newline 
3) We first note that via \eqref{E:homo} the $G$-action on $\Hom(\Gamma, G)$ 
is compatible with the diagonal $G$-action on $\tG$: 
$$(g_1, \cdots, g_k) \; 
\mapsto \; (h g_1 h^{-1}, \cdots, h g_k h^{-1}).$$ 
Now, we compute the $j$th components of the 
image of $\Psi_1$ (see \eqref{E:e1} for the definition) and 
$\Psi_2$ (see \eqref{E:e2} for the definition), respectively, 
as follows: 
\newline
For $h = \begin{pmatrix} 
I_k & \vec{a} & \vec{b}\\
0 & 1 & c \\
0 & 0 & 1
\end{pmatrix} \in G$, we have 
$$
h \exp \begin{pmatrix} \bigzero  & z_j \ya{x} & \ya{y_j} \\
              0 & 0 & z_j \\
              0 & 0 & 0
     \end{pmatrix} h^{-1} = 
\exp \begin{pmatrix} \bigzero  & z_j \ya{x} & \ya{y_j} + z_j (\ya{a} - c \ya{x}) \\
              0 & 0 & z_j \\
              0 & 0 & 0
     \end{pmatrix} 
$$ 
and 
$$
h \exp \begin{pmatrix} \bigzero  & \ya{x_j} & \ya{y_j} \\
              0 & 0 & 0 \\
              0 & 0 & 0
     \end{pmatrix} h^{-1} = 
\exp \begin{pmatrix} \bigzero  & \ya{x_j} & \ya{y_j} -c \ya{x_j} \\
              0 & 0 & 0 \\
              0 & 0 & 0
     \end{pmatrix} 
$$ 
This is what we wanted to prove. 
\newline
4) This is an elementary linear algebra.
\end{proof}

Then, Theorem \ref{T:deformspace} 
is an immediate consequence of
 Theorem~\ref{T:deformpara} and Proposition~\ref{P:keyprop} (4).
Corollary~\ref{cor:dimdefo} now follows from
 Proposition~\ref{P:2.4.1},
and from the $G$-action on $R(\Gamma,G;\mathbb{R}^{k+1})$
described in the above proof.

\bigskip

\textbf{Acknowledgement:} \ 
Part of the results here was obtained
 while the first author was invited to
Harvard University 2000--2001 and the second author was in the  
University of Tokyo,
and  the final version was made
 while both of them were working at
 RIMS, Kyoto University.
We would like to thank these institutions for warm atmosphere of research.
A preliminary version of this paper was circulated as
``Proper action of ${\mathbb R}^k$ on a $(k+1)$-dimensional
nilpotent homogeneous manifold".


\begin{thebibliography}{99}
\bibitem{xams}
{\sc H. Abels, G. A. Margulis and G. A. Soifer},
{\it On the Zariski closure of the linear part of a properly discontinuous group of affine transformations},
J. Differential Geom. {\bf 60} (2002), 315--344.

\bibitem{baklouti} 
{\sc A. Baklouti and F. Khlif}, 
{\it Proper actions on some exponential solvable homogeneous spaces}, 
Internat. J. Math. {\bf 16} (2005), 941--955.

\bibitem{benoist}
{\sc Y. Benoist}, 
{\it Actions propres sur les espaces homog\`{e}nes r\'{e}ductifs}, 
Ann. of Math. {\bf 144} (1996), 315--347.
 
 
\bibitem{goldman}
{\sc W. Goldman}, 
{\it Nonstandard Lorentz space forms}, 
J. Differential Geom. {\bf 21} (1985), 301--308.


\bibitem{kobayashi89}
{\sc T. Kobayashi}, 
{\it Proper action on a homogeneous space of reductive type}, 
Math. Ann. {\bf 285} (1989), 249--263.


\bibitem{kobayashi93}
\bysame, 
{\it On discontinuous groups acting on homogeneous spaces with noncompact isotropy subgroups}, 
J. Geom. Phys. {\bf 12} (1993), 133--144.

\bibitem{kobayashi96}
\bysame, 
{\it Criterion for proper actions on homogeneous spaces of reductive groups}, 
J. Lie Theory {\bf 6} (1996), 147--163.


\bibitem{kobayashi98}
\bysame, 
{\it Deformation of compact Clifford--Klein forms of indefinite-Riemannian homogeneous manifolds}, 
Math. Ann. {\bf 310} (1998), 395--409.

\bibitem{kobayashi00}
\bysame, 
{\it Discontinuous groups for non-Riemannian homogeneous spaces}, 
Mathematics Unlimited --- 2001 and Beyond 
(B. Engquist and W. Schmid, eds.), Springer, 2001, 723--747.

\bibitem{kobayashi05}
\bysame,
{\it On discontinuous group actions on non-Riemannian homogeneous spaces},
 Sugaku {\bf 57} (2005), 267--281
  (in Japanese);
 An English translation is to appear in 
 ``Sugaku Expositions'', Amer. Math. Soc., 
preprint RIMS-1537.
 

\bibitem{lipsman}
{\sc R. Lipsman}, 
{\it Proper actions and a compactness condition}, 
J. Lie Theory {\bf 5} (1995), 25--39.

\bibitem{nasrin1}
{\sc S. Nasrin}, 
{\it On a conjecture of Lipsman about proper actions
on nilpotent Lie groups}, 
Master thesis, the University of Tokyo, January 2000.

\bibitem{nasrin}
\bysame,   % {\sc Nasrin, ~S.} 
{\it Criterion of proper actions for 2-step nilpotent Lie groups}, 
Tokyo J. Math. {\bf 24} (2001), 535--543.

\bibitem{palais}
{\sc R. S. Palais}, 
{\it On the existence of slices for actions of non-compact Lie groups}, 
Ann. of Math. {\bf 73} (1961), 295--323.


\bibitem{salein}
{\sc F. Salein}, 
{\it Vari\'{e}t\'{e}s anti-de Sitter de dimension 3 poss\'{e}dant un champ de Killing non trivial}, 
C. R. Acad. Sci. Paris {\bf 324} (1997), 525--530; 
Th\`{e}se \`{a} l'\'{E}cole Normale Sup\'{e}rieure de Lyon, December 1999.

\bibitem{weil}
{\sc A. Weil}, 
{\it Remarks on the cohomology of groups}, 
Ann. of Math. {\bf 80} (1964), 149--157.


\bibitem{witte}
{\sc D. Witte},
{\it Superrigid subgroups and syndetic hulls in solvable Lie groups},
Rigidity in dynamics and geometry (M. Burger and A. Iozzi, eds.)
(Cambridge, 2000),
Springer, 2002, 441--457.

\bibitem{yoshino3step}
{\sc T. Yoshino}, 
{\it Criterion of proper actions for three step nilpotent Lie groups}, 
preprint.

\bibitem{yoshino}
\bysame, %{\sc Yoshino,~T.},
 personal communication.
 
\end{thebibliography}
\end{document}